\documentclass{article}

\usepackage{amsmath}
\usepackage{amssymb}
\usepackage{amsfonts}
\usepackage{amsthm}
\usepackage{cleveref}
\usepackage{fullpage}
\usepackage{xcolor}
\usepackage{graphicx}

\newcommand{\R}{\mathbb{R}}

\renewcommand{\d}{\mathrm{d}}
\renewcommand{\div}{\mathrm{div}}

\newcommand{\Poisson}{\mathrm{P}}
\newcommand{\EEG}{\mathrm{EEG}}

\newtheorem{theorem}{Theorem}
\newtheorem{proposition}[theorem]{Proposition}

\graphicspath{{./figures/},{.}}

\title{Green's function for Poisson's equation and the EEG equation with Neumann boundary condition on $n$-balls}
\author{Benedikt Wirth}
\date{}

\begin{document}

\maketitle

\begin{abstract}
We provide an elementary derivation of the Green's function for Poisson's equation with Neumann boundary data on balls of arbitrary dimension, which was recently found in \cite{SaToTu16}.
The underlying idea consists of first computing the Green's function for the electroencephalography (EEG) equation (Poisson's equation with dipole right-hand side)
and then deriving the Green's function for Poisson's equation from that.
\end{abstract}

\section{Introduction}
A classic approach to solve linear partial differential equations is by the use of Green's functions.
In simple geometries, such as halfplanes and balls, and for simple partial differential equations these can often be stated explicitly.
Surprisingly, for the Neumann--Poisson problem in balls of arbitrary dimension, the Green's function was only derived recently \cite{SaToTu16} (it had previously been known only in dimensions 1, 2, and 3).
Unaware of this result, the author of this note studied Poisson's problem with a dipole right-hand side and Neumann boundary data (the so-called forward problem of electroencephalography)
and constructed a Green's function for this problem by standard methods.
The intimate connection to the Neumann--Poisson problem then provides a direct way to find its Green's function as well.

In principle, once a Green's function has been stated, it can be verified by straightforward calculation.
However, how the Green's function was found may be of interest as well, be it as a source of inspiration for similar constructions or simply to give some insight into an otherwise magic formula.
For the Neumann--Poisson problem, one of the difficulties is that the boundary conditions have to be compatible with the right-hand side of Poisson's equation.
The derivation given below circumvents this difficulty by first studying a dipole right-hand side, for which it suffices to apply the standard reflection method. 
In contrast, the derivation in \cite{SaToTu16} (which follows the three-dimensional construction from \cite[p.\,348]{BoGiDe35}) uses infinite series expansions with harmonic polynomials of the fundamental solution for Poisson's equation, which are more complicated and less intuitive.

In summary, below we will find (distributional) solutions $G_z^\Poisson$ and $G_z^\EEG$ (the sought Green's functions) to the Poisson--Neumann problem
\begin{align}
\Delta G_z^\Poisson&=\delta_z-\frac{n}{\omega_n}&&\text{in }B\,,\label{eqn:Poisson}\\
\nabla G_z^\Poisson\cdot\nu&=0&&\text{on }\partial B\label{eqn:PoissonBC}
\end{align}
and to the EEG problem
\begin{align}
\Delta G_z^{\EEG,D}&=\div(D\delta_z)&&\text{in }B\,,\label{eqn:EEG}\\
\nabla G_z^{\EEG,D}\cdot\nu&=0&&\text{on }\partial B\,,\label{eqn:EEGBC}
\end{align}
where $B$ denotes the unit ball in $\R^n$ and $\omega_n$ its surface area, $\nu\in\R^n$ denotes the unit outward normal on $\partial B$,
$\delta_z$ is the Dirac distribution at $z\in B$, $D\in\R^n$ is the direction of the dipole, and $\div$ represents the (distributional) divergence.
The constant term $-\frac n{\omega_n}$ on the right-hand side of \eqref{eqn:Poisson} represents the inverse unit ball volume and is necessary to ensure the solvability of the Neumann--Poisson problem.
The corresponding Green's functions $G_z^{\Poisson,R}$ and $G_z^{\EEG,D,R}$ on balls of radius $R\neq1$ are then simply obtained as
\begin{align*}
G_z^{\Poisson,R}(x)&=R^{n-2}G_{z/R}^\Poisson(x/R)\,,\\
G_z^{\EEG,D,R}(x)&=R^{n-1}G_{z/R}^{\EEG,D}(x/R)\,,
\end{align*}
as follows from straightforward scaling arguments.

\section{Green's function for the EEG equation with radial dipoles}
Let $e$ be a unit vector in $\R^n$ and $c\in[0,1)$.
In this section we provide $G_z^{\EEG,D}$ for $z=ce$ and $D=e$, that is, a dipole at $z$ pointing radially outwards.
To this end, recall the fundamental solution for the Laplace operator
\begin{equation*}
\Phi(x)=\begin{cases}
\frac1{\omega_1}|x|&\text{if }n=1,\\
\frac1{\omega_2}\log|x|&\text{if }n=2,\\
-\frac1{(n-2)\omega_n|x|^{n-2}}&\text{if }n>2,
\end{cases}
\end{equation*}
which satisfies $\Delta\Phi=\delta_0$.
By taking the gradient on both sides we obtain $\Delta\Psi=\nabla\delta_0$ for the fundamental solution
\begin{equation*}
\Psi(x)=\nabla\Phi(x)=\frac{x}{\omega_n|x|^n}
\end{equation*}
of the EEG problem.
Thus, abbreviating $\Psi^D=D\cdot\Psi$ we have
\begin{equation*}
\Delta\Psi^D(\cdot-z)=D\cdot\nabla\delta_z=\div(D\delta_z)\,,
\end{equation*}
however, $\Psi^D$ does not yet fulfill the homogeneous Neumann boundary conditions \eqref{eqn:EEGBC}.
Those are enforced by the well-known reflection method.
To this end let $z^*=z/|z|^2$ be the inverse of $z$ with respect to the sphere $\partial B$ and define
\begin{equation*}
G_z^{\EEG,D}(x)
=\Psi^D(x-z)+\Psi^{-D}(x-z^*)/|z|^n
=\frac D{\omega_n}\cdot\left(\frac{x-z}{|x-z|^n}-\frac{x-z^*}{|z|^n|x-z^*|^n}\right)\,,
\end{equation*}
then indeed for $D=e$ and $z=ce$ we have $\Delta G_z^{\EEG,D}=\div(D\delta_z)$ in $B$ and
\begin{equation*}
\nabla G_z^{\EEG,D}(x)\cdot\nu
=\frac e{\omega_n}\left(\frac1{|x-z|^n}-\frac1{|z|^n|x-z^*|^n}\right)\cdot\nu-\frac n{\omega_n}\left(\frac{e\cdot(x-z)}{|x-z|^{n+2}}(x-z)-\frac{e\cdot(x-z^*)}{|z|^n|x-z^*|^{n+2}}(x-z^*)\right)\cdot\nu
=0
\end{equation*}
on $\partial B$, where we exploited $\nu=x$ on $x\in\partial B$ and the relation $|x-z|=|z||x-z^*|$ central to the reflection method.
In summary,
\begin{equation*}
G_{ce}^{\EEG,e}(x)
=\frac1{\omega_n}\left(\frac{e\cdot x-c}{|x-ce|^n}-\frac{e\cdot x-\frac1c}{c^n|x-\frac ec|^n}\right)\,.
\end{equation*}

\section{Green's function for the Neumann--Poisson problem}
Upon differentiating \eqref{eqn:Poisson}-\eqref{eqn:PoissonBC} with respect to $z$ in direction $-D\in\R^n$, we obtain
\begin{align*}
\Delta(-D\cdot\nabla_zG_z^\Poisson)&=\div(D\delta_z)&&\text{in }B\,,\\
\nabla(-D\cdot\nabla_zG_z^\Poisson)\cdot\nu&=0&&\text{on }\partial B\,,
\end{align*}
where $\nabla_z$ shall denote the gradient with respect to $z\in\R^n$.
Consequently we have $G_z^{\EEG,D}=-D\cdot\nabla_zG_z^\Poisson$, which is the relation from which in the next section we will obtain the Green's function for the EEG equation with arbitrary dipoles.
Here, we just exploit the special case
\begin{equation*}
G_{ce}^{\EEG,e}=-\frac\partial{\partial c}G_{ce}^\Poisson\,,
\end{equation*}
which is obtained by choosing $z=ce$ and $D=e$ in the above for $e\in\R^n$ a unit vector and $c\in[0,1)$.
It shows that $G_{ce}^\Poisson$ is an antiderivative of $-G_{ce}^{\EEG,e}$, which we computed in the previous section.
Hence, in order to find $G_{ce}^\Poisson$ it is natural to introduce the auxiliary function
\begin{equation*}
F_c(x)=-\int_0^cG_{se}^{\EEG,e}(x)-\frac1{\omega_ns}\,\d s\,.
\end{equation*}
Here we exploited that an arbitrary function of $s$ can be added to $G_{se}^{\EEG,e}$ without violating its Green's function property; the choice $\frac1{\omega_ns}$ just makes the integral finite.
By construction, the function $F_c$ satisfies
\begin{align*}
\Delta F_c
&=-\int_0^c\Delta G_{se}^{\EEG,e}\,\d s
=-\int_0^c\div(e\delta_{se})\,\d s
=\int_0^c\frac\partial{\partial s}\delta_{se}\,\d s
=\delta_{ce}-\delta_0
&&\text{in }B\,,\\
\nabla F_c\cdot\nu
&=-\int_0^c\nabla G_{se}^{\EEG,e}\cdot\nu\,\d s
=0
&&\text{on }\partial B\,.
\end{align*}
Furthermore, it is straightforward to calculate the radially symmetric solution
\begin{equation*}
G_0^\Poisson(x)=\Phi(x)-\frac1{2\omega_n}|x|^2
\end{equation*}
to \eqref{eqn:Poisson}-\eqref{eqn:PoissonBC}
(indeed, under radial symmetry, \eqref{eqn:Poisson}-\eqref{eqn:PoissonBC} turn into a readily solvable ordinary differential boundary value problem).
As a consequence, we have
\begin{multline*}
G_{ce}^\Poisson(x)
=F_c(x)+G_0^\Poisson(x)
=\frac{-1}{\omega_n}\int_0^c\frac{e\cdot x-s}{|x-se|^n}-\frac{e\cdot x-\frac1s}{s^n|x-\frac es|^n}-\frac1s\,\d s
+\Phi(x)-\frac1{2\omega_n}|x|^2\\
=\Phi(x-ce)
+\frac1{\omega_n}\int_0^c\frac{e\cdot x-\frac1s}{s^n|x-\frac es|^n}+\frac1s\,\d s
-\frac1{2\omega_n}|x|^2\,.
\end{multline*}
This is the desired formula for the Green's function, thus we can summarize as follows.

\begin{proposition}[Green's function for Neumann--Poisson problem on the unit ball]
Equations \eqref{eqn:Poisson}-\eqref{eqn:PoissonBC} are solved up to an arbitrary additive function of $z$ by
\begin{align*}
G_z^\Poisson(x)
&=\Phi(x-z)
+\frac1{\omega_n}\int_0^{|z|}\frac{x\cdot\frac z{|z|}-\frac1s}{|sx-\frac z{|z|}|^n}+\frac1s\,\d s
-\frac1{2\omega_n}|x|^2\,.
\end{align*}
\end{proposition}

The integral within the formula can be analytically calculated using just elementary functions; for the sake of completeness and for reference we provide the corresponding (lengthy) formulas in the appendix.

\section{Green's function for the EEG equation}
Having computed $G_z^\Poisson$, we obtain the Green's function for the EEG equation via $G_z^{\EEG,D}=-D\cdot\nabla_zG_z^\Poisson$, as discussed in the previous section.

\begin{proposition}[Green's function for EEG problem on the unit ball]
Equations \eqref{eqn:EEG}-\eqref{eqn:EEGBC} are solved up to an arbitrary additive function of $z$ by
\begin{align*}
G_z^{\EEG,D}(x)
&=\frac D{\omega_n}\cdot\left(\frac{(x-z)}{|x-z|^n}
-\frac{x\cdot z-1}{||z|x-\frac z{|z|}|^n}\frac z{|z|^2}-\frac z{|z|^2}
-\int_0^{|z|}\frac1{|sx-\frac z{|z|}|^n}+n\frac{sx\cdot\frac z{|z|}-1}{|sx-\frac z{|z|}|^{n+2}}\,\d s\left(I-\frac z{|z|}\otimes\frac z{|z|}\right)\frac x{|z|}\right)\\
&=\frac D{\omega_n}\cdot\Bigg(\frac{(x-z)}{|x-z|^n}
-\frac{x\cdot z-1}{|z|^n|x-z^*|^n}z^*-z^*
-\Bigg[\frac{x\cdot\frac z{|z|}}{|x|^2}\Bigg(\frac1{|z|^n|x-z^*|^n}-1\Bigg) \\
&\qquad+\int_0^{|z|}\frac1{|sx-\frac z{|z|}|^n}\,\d s+n\left(\left(\frac x{|x|}\cdot\frac z{|z|}\right)^2-1\right)\int_0^{|z|}\frac{1}{|sx-\frac z{|z|}|^{n+2}}\,\d s\Bigg]\left(I-\frac z{|z|}\otimes\frac z{|z|}\right)\frac x{|z|}\Bigg)\,.
\end{align*}
\end{proposition}

Again the integrals can be computed analytically and expressed via elementary functions; the corresponding formulas are provided in the appendix.

\appendix
\section{Appendix}
The formula for $G_z^\Poisson$ requires the evaluation of the integral $\int_0^{|z|}\frac{x\cdot z/|z|-\frac1s}{|sx-\frac z{|z|}|^k}+\frac1s\,\d s$ for positive integers $k$.
Abbreviating $e=z/|z|$ and $c=|z|$, the integral to be computed reads
\begin{equation*}
\Gamma_k(x,e,c)=\int_0^c\frac{e\cdot x-\frac1s}{|sx-e|^k}+\frac1s\,\d s\,.
\end{equation*}
We have the following expressions (below, $k>2$),
\begin{align*}
\Gamma_1(x,e,c)
&=\log\left(\frac c2\sqrt{|x|^2-(e\cdot x)^2}\right)+\frac{e\cdot x}{|x|}\log\left(\frac{c|x|^2-e\cdot x+|x||cx-e|}{|x|-e\cdot x}\right)-\mathrm{arctanh}\left(\frac{ce\cdot x-1}{|cx-e|}\right)\,,\\
\Gamma_2(x,e,c)
&=\log|cx-e|\,,\\
\Gamma_k(x,e,c)
&=\left[\frac{-1}{(k-2)|sx-e|^{k-2}}\right]_{s=0}^c-\int_0^c\frac{\frac1s}{|sx-e|^{k-2}}-\frac1s\,\d s\\
&=\frac1{k-2}-\frac1{(k-2)|cx-e|^{k-2}}+\Gamma_{k-2}(x,e,c)-\int_0^c\frac{e\cdot x}{|sx-e|^{k-2}}\,\d s\\
&=\sum_{\substack{j=k-2,k-4,\ldots\\j>0}}\left(\frac1j-\frac1{j|cx-e|^j}-e\cdot x\int_0^c\frac1{|sx-e|^j}\,\d s\right)+\begin{cases}\Gamma_2(x,e,c)&\text{if $k$ even,}\\\Gamma_1(x,e,c)&\text{if $k$ odd,}\end{cases}
\end{align*}
where the last line follows from resolving the recursion.
Note that in one spatial dimension, $\Gamma_1(x,e,c)$ simplifies to zero.
It remains to calculate
\begin{equation*}
Z_j
=\int_0^c\frac1{|sx-e|^j}\,\d s
=\int_0^c\frac1{\sqrt{s^2|x|^2-2se\cdot x+1}^j}\,\d s
=\int_0^c\frac1{\sqrt{(s|x|-\frac{e\cdot x}{|x|})^2+1-(\frac{e\cdot x}{|x|})^2}^j}\,\d s\,.
\end{equation*}
Abbreviating $A^2=1-(\frac{e\cdot x}{|x|})^2$ and changing variables to $t=(s|x|-\frac{e\cdot x}{|x|})/A$ we obtain
\begin{equation*}
Z_j
=\frac1{|x|A^{j-1}}\int_{-\frac{e\cdot x}{A|x|}}^{\frac{c|x|}A-\frac{e\cdot x}{A|x|}}\frac1{\sqrt{t^2+1}^j}\,\d t
=\frac1{|x|A^{j-1}}\left(\int_{0}^{\frac{c|x|}A-\frac{e\cdot x}{A|x|}}\frac1{\sqrt{t^2+1}^j}\,\d t-\int_0^{-\frac{e\cdot x}{A|x|}}\frac1{\sqrt{t^2+1}^j}\,\d t\right)
\,.
\end{equation*}
We finish by calculating $\int_0^a(1+t^2)^{-j/2}\,\d t$, where without loss of generality we assume $a\geq0$ (otherwise we may replace $a$ by $|a|$ and thereby only change the sign of the integral).
We first consider the case of odd $j$.
By substituting $t=1/\sqrt{z^2-1}$ such that $\frac{\d t}{\d z}=-\frac{z}{\sqrt{z^2-1}^3}$ we obtain
\begin{equation*}
\int_0^a\frac1{\sqrt{1+t^2}^j}\,\d t
=\int_{1+1/a^2}^\infty z^{1-j}\sqrt{z^2-1}^{j-3}\,\d z\,.
\end{equation*}
For $j=1$ we thus obtain
\begin{equation*}
\int_0^a\frac1{\sqrt{1+t^2}^j}\,\d t
=\int_{1+1/a^2}^\infty\frac1{z^2-1}\,\d z
=\left[\frac12\log\frac{z-1}{z+1}\right]_{z=1+1/a^2}^\infty
=\frac12\log(2a^2+1)
\,,
\end{equation*}
and for odd $j>1$ we have
\begin{multline*}
\int_0^a\frac1{\sqrt{1+t^2}^j}\,\d t
=\int_{1+1/a^2}^\infty\frac{(z^2-1)^{\frac{j-3}2}}{z^{j-1}}\,\d z\\
=\int_{1+1/a^2}^\infty\sum_{i=0}^{\frac{j-3}2}(-1)^{i}{\frac{j-3}2\choose i}z^{2i+1-j}\,\d z
=\sum_{i=0}^{\frac{j-3}2}\frac{(-1)^i}{2i+2-j}{\frac{j-3}2\choose i}\left(1+\frac1{a^2}\right)^{2i+2-j}
\,.
\end{multline*}
Now we provide the case of even $j$,
\begin{align*}
\int_0^a\frac1{\sqrt{1+t^2}^j}\,\d t
&=\left[\frac t{(j-2)\sqrt{1+t^2}^{j-2}}\right]_{t=0}^a+\frac{j-3}{j-2}\int_0^a\frac1{\sqrt{1+t^2}^{j-2}}\,\d t\\
&=\frac{(j-3)!!}{(j-2)!!}\left(\arctan(a)+\sum_{i=j-2,j-4,\ldots,2}\frac{i!!}{(i-1)!!}\frac a{i\sqrt{1+a^2}^i}\right)\,,
\end{align*}
where in the second step we resolved the recursion and $k!!=k(k-2)(k-4)\cdots$ denotes the double factorial.

\paragraph{Acknowledgements.}
This work has been supported by by the Alfried Krupp Prize for Young University Teachers awarded by the Alfried Krupp von Bohlen und Halbach-Stiftung
as well as by the Deutsche Forschungsgemeinschaft (DFG, German Research Foundation) under Germany's Excellence Strategy -- EXC 2044 --, Mathematics M\"unster: Dynamics -- Geometry -- Structure.

\bibliographystyle{plain}
\bibliography{notes}

\end{document}